\documentclass{amsart}
\usepackage{amssymb,latexsym,graphics,amsfonts}

\swapnumbers
 \theoremstyle{plain}

 \theoremstyle{definition}
 
 \theoremstyle{remark}
 
 \newcommand{\A}{\mathcal{A}}

 \newcommand{\cal}[1]{\mathcal{#1}}

\begin{document}
\title{Module extensions of dual Banach algebras}
\author{M. Eshaghi Gordji}
\address{Department of Mathematics,
University of Semnan, Semnan, Iran} \email{maj\_ess@Yahoo.com and
meshaghi@semnan.ac.ir}
\author{ F. Habibian}
\address{Department of Mathematics,
Isfahan University, Isfahan, Iran}
\email{fhabibian@ui.ac.ir}
\author{A. Rejali}
\address{Department of Mathematics,
Isfahan University, Isfahan, Iran}\email{rejali@ui.ac.ir}
\dedicatory{}
\smallskip
\begin{abstract}
In this paper we define the module extension dual Banach algebras
and we use this Banach algebras to finding the relationship between
$weak^*-$continuous homomorphisms of dual Banach algebras and
Connes-amenability. So we study the $weak^*-$continuous derivations
on module extension dual Banach algebras.
\end{abstract}
\maketitle
$$\bf Introduction$$
A Banach algebra $\cal A$ is said to be dual if there is a closed
submodule ${\cal A}_*$ of $\cal A^*$ such that $\cal A={{\cal
A}_*}^*.$ Let $\cal A$ be a dual Banach algebra. A dual Banach $\cal
A$-bimodule $X$ is called normal if, for every $x\in X,$ the maps
$a\longmapsto a.x$ and
$a\longmapsto x.a$ are $weak^*-$continuous from $\cal A$ into $X$. \\
For example if $G$ is a locally compact topological group, then
$M(G)$ is a dual Banach algebra with predual $C_0(G)$. Also if $\cal
A$ is an Arens regular Banach algebra, then $\cal A^{**}$ (by the
first (or second) Arens product) is a dual Banach algebra with
predual $\cal A^*.$
 Let $\cal A$ and $\cal B$ be dual Banach algebras and let $\varphi:\cal
A\longrightarrow\cal B$ be a $weak^*-$continuous Banach algebra
homomorphism, then $\cal B$ is a normal $\cal A$-bimodule by the
following module actions
$$a.b=\varphi(a)b , \hspace{0.7cm} b.a=b\varphi(a) \hspace{1.5cm} (a\in
\cal A, b\in \cal B).$$ We denote $\cal B_{\varphi}$ the above
$\cal A$-bimodule. If $X$ is a Banach $\cal A$-bimodule then a
derivation from $\cal A$ into $X$ is a linear map $D$, such that
for every $a,b \in \cal A,$ $D(ab)=D(a).b+a.D(b).$ If $x\in X,$
and we define $\delta_x:\cal A\longrightarrow X$ by
$\delta_x(a)=a.x-x.a\hspace {0.5cm}(a\in \cal A),$ then $\delta_x$
is a derivation, derivations of this form are called inner
derivations. Let $\cal A$ be a Banach algebra and $X$ be a Banach
$\cal A$-module, then $X^*$ is a Banach $\cal A$-module if for
each $a\in \cal A$ and $x\in X$ and $x^*\in X^*$ we define
\[\langle ax^*,x \rangle=\langle x^*,xa \rangle, \hspace {1 cm} \langle x^*a,x \rangle=\langle x^*,ax\rangle.\]
A Banach algebra $\cal A$  is amenable if every derivation from
$\cal A$ into every dual $\cal A$-module is inner, equivalently if
$H^1(\cal A,X^*)=\{0\}$ for every Banach $A$-module $X$, where
$H^1(\cal A,X^*)$ is the first cohomology group from $\cal A$ with
coefficients in $X^*$[3] (see [1] and [5], for more details).\\A
dual Banach algebra $\cal A$ is Connes-amenable if every
$weak^*-$continuous derivation from $\cal A$ into each normal dual
Banach $\cal A$-bimodule $X$
 is inner; i.e. $H^1_{w^*}(\cal A, X)=\{o\}$, this definition was introduced by V. Runde
(see section 4 of [5]). A dual Banach algebra $\cal A$ is weakly
Connes-amenable if every $weak^*-$continuous derivation from $\cal
A$ into $\cal A$ is inner; i.e. $H^1_{w^*}(\cal A, \cal A)=\{o\}$
[2]. Yong Zhang studied the weak amenability of module extension
Banach algebras [6]. We define the module extensions of dual Banach
algebras and then we study the Connes-amenability and the weak
Connes-amenability of Banach algebras of this form.

\section{Connes  Amenability}
In this section we  find a necessary and sufficient condition for a
dual Banach algebra to be Connes-amenable. First we define a new
class of dual Banach algebras. Let $\cal A$ be a dual Banach algebra
with predual $\cal A_*$, and let $X$ be a normal dual Banach $\cal
A$-bimodule by predual $X_*$. Then we can show that ${\cal
A}\oplus_{\infty} X$  is a Banach space with the following norm
$$\|(a,x)\|=Max\{\|x\|,\|a\|\}\hspace{1.5cm} (a\in{\cal A}~,~x\in X).$$
So $\cal A \oplus_ {\infty}X$ is a Banach algebra with the following
product ,
$$(a_1,x_1)(a_2,x_2)=(a_1a_2,a_1\cdot x_2+x_1\cdot a_2)~ .$$
${\cal A}_*\oplus_1 X_*$  is a Banach space with the norm
$\|(a',x')\|=\|a'\|+\|x'\| \hspace{0.5cm} (a'\in{\cal A_*}~,~x'\in
X_*).$ We have  ${\cal A}\oplus_{\infty}X =({{\cal A}_*\oplus_1
X_*})^*$. Since $\cal A$ is a dual Banach algebra and $X$ is a
normal dual Banach $\cal A-$bimodule, then it is easy to show that
the multiplication in $\cal A \oplus_ {\infty}X$ is separately
$weak^*-$continuous. Thus by 4.4.1 of [5], we have the following
lemma.
\paragraph{\bf Lemma 1.1.} Let $\cal A$ and $X$ be as above, then
$\cal A \oplus_ {\infty}X$ is a dual Banach algebra with predual
${\cal A}_*\oplus_ {1}X_*.$

The Banach algebra $\cal A \oplus_ {\infty}X$ in lemma is called
module extension dual Banach algebra.
\paragraph{\bf Theorem 1.2.}
Let $\cal A$ be a dual Banach algebra. Then the following assertions
are equivalent:

(i) $\cal A$ is Connes-amenable.

(ii) For every dual  Banach algebra $\cal B$ and every
$weak^*-$continuous homomorphism $\varphi:\cal A\longrightarrow\cal
B$, $H^1_{w^*}(\cal A,\cal {B_{\varphi}})=\{o\}.$

(iii) For every dual  Banach algebra $\cal B$ and every injective
$weak^*-$continuous homomorphism $\varphi:\cal
A\longrightarrow\cal B$, $H^1_{w^*}(\cal A,\cal
{B_{\varphi}})=\{o\}.$
\paragraph{\bf Proof.} The proofs of $(i)\Longrightarrow (ii)$ and $(ii)\Longrightarrow (iii)$ are
easy. We prove $(iii)\Longrightarrow (i)$. Let $X$ be a normal dual
Banach $\cal A-$bimodule, and let $D:\cal A\longrightarrow X$ be a
$weak^*-$continuous derivation. Obviously by above lemma, the map
$$\varphi:  a\mapsto(a,0),\hspace{0.7cm} \cal A\longrightarrow
\cal A\oplus_{\infty} X$$ is an injective $weak^*-$continuous
homomorphism. Thus $H^1_{w^*}(\cal A, ({(\cal A \oplus_{\infty}
X)_{\varphi}}))=\{o\}.$ We define $D_1:\cal A\longrightarrow \cal A
\oplus_{\infty} X$ by $D_1(a)=(0,D(a))$. For $a,b \in \cal A$ we
have
\begin{align*} D_1(ab)&=(0,D(ab))=(0,D(a)b+aD(b)) \\
&=(0,D(a))(b,0)+(a,0)(0,D(b))\\
&=D_1(a)\varphi(b)+\varphi(a)D_1(b).
\end{align*}
Thus $D_1$ is a $weak ^*-$continuous derivation from $\cal A$ into
${(\cal A \oplus_{\infty} X)_{\varphi}}$, and then it is inner
derivation. On the other word there exist $b\in \cal A , x\in X$
such that $D_1=\delta_{(b,x)}.$ For every $a\in \cal A$ we have
\begin{align*} (0,D(a))&=D_1(a)=\delta_{(b,x)}(a) \\
&=\varphi(a)(b,x)-(b,x)\varphi(a) \\
&= (a,0)(b,x)-(b,x)(a,0) \\
&=(ab-ba,ax-xa).
\end{align*}
Thus $D=\delta_{x}$. So $\cal A$ is
Connes-amenable.\hfill$\blacksquare~$

Now we find the necessary and sufficient condition for a dual module
extension Banach algebra to be Connes-amenable.
\paragraph{\bf Theorem 1.3.}
Let $\cal A$ be a dual Banach algebra and let $X$ be a reflexive
Banach $\cal A$-bimodule. If for every $x'\in X^*$ and $a\in \cal
A,$ the mappings
$$(x'\widehat{\otimes}a).:b\longmapsto (x'\widehat{\otimes}ab)~ ,\hspace{0.5cm}
 .(x'\widehat{\otimes}a):b\longmapsto (bx'\widehat{\otimes}a); \hspace {0.5cm} \cal A
\longrightarrow X^*\widehat{\otimes} {\cal A},\hspace {0.5cm}(1)$$
are $weak^*-$weak continuous, then $\cal A \oplus_ {\infty}X$  is
Connes-amenable if and only if $X=0$ and $\cal A$ is
Connes-amenable.
\paragraph{\bf Proof.} Let  $\cal A \oplus_ {\infty}X$  be Connes-amenable and the mappings
defined in (1), are $weak^*-$weak continuous, we have to show that
$X=0$. It is easy to check that $X^*\widehat{\otimes} {\cal A}$ is a
Banach $\cal A \oplus_ {\infty}X$-bimodule with the following module
actions:
$$(x'\widehat{\otimes} a).(b,x)=x'\widehat{\otimes}ab, \hspace{0.4cm}(b,x)
.(x'\widehat{\otimes}a)=bx'\widehat{\otimes}a,\hspace{0.7cm}
(x'\widehat{\otimes}a \in {X^*\widehat{\otimes} {\cal A}}, (b,x)\in
\cal A \oplus_ {\infty}X).
$$
 Let
$$(b_\alpha,x_\alpha)~\stackrel{weak^*}{-\hspace{-.2cm}-\hspace{-.2cm}\longrightarrow}
(b,x) \hspace{1cm} in \hspace{0.3cm} \cal A \oplus_ {\infty}X,$$
thus
$b_\alpha~\stackrel{weak^*}{-\hspace{-.2cm}-\hspace{-.2cm}\longrightarrow}
b$ in $\cal A$. Then for each $x'\in X^*$ and each $a\in \cal A,$ we
have $$b_\alpha
x'\widehat{\otimes}a~\stackrel{weakly}{-\hspace{-.2cm}-\hspace{-.2cm}\longrightarrow}
bx'\widehat{\otimes}a \hspace{1cm} in \hspace{0.3cm}
X^*\widehat{\otimes} {\cal A}.$$ For each  $F\in
{(X^*\widehat{\otimes} {\cal A})}^*$, we have
$$\langle F.(b_\alpha,x_\alpha),x'\widehat{\otimes}a\rangle= \langle
F,b_\alpha x'\widehat{\otimes}a\rangle \longrightarrow \langle F,
bx'\widehat{\otimes}a \rangle=\langle
F.(b,x),x'\widehat{\otimes}a\rangle.$$ This means that
$$F.(b_\alpha,x_\alpha)~\stackrel{weak^*}{-\hspace{-.2cm}-\hspace{-.2cm}\longrightarrow}
F.(b,x)\hspace{1cm} in \hspace{0.3cm} {(X^*\widehat{\otimes} {\cal
A})}^*.$$ Similarly we have
$$(b_\alpha,x_\alpha).F~\stackrel{weak^*}{-\hspace{-.2cm}-\hspace{-.2cm}\longrightarrow}
(b,x).F \hspace{1cm} in \hspace{0.3cm} {(X^*\widehat{\otimes} {\cal
A})}^*.$$ Thus ${(X^*\widehat{\otimes} {\cal A})}^*$ is a normal
dual $\cal A \oplus_ {\infty}X$-bimodule. We define $D:\cal A
\oplus_ {\infty}X \longrightarrow {(X^*\widehat{\otimes} {\cal
A})}^*$ as follows:
$$\langle D(b,x),x'\widehat{\otimes}a\rangle=\langle x',ax\rangle \hspace{0.7cm} (x'\widehat{\otimes}a \in
{X^*\widehat{\otimes} {\cal A}}, (b,x)\in \cal A \oplus_
{\infty}X).$$ For each $(b_1,x_1), (b_2,x_2)\in \cal A \oplus_
{\infty}X,$ and $x'\widehat{\otimes}a \in {X^*\widehat{\otimes}
{\cal A}},$ we have
\begin{align*} \langle D((b_1,x_1)(b_2,x_2)),x'\widehat{\otimes}a \rangle&=\langle x',a(b_1x_2+x_1b_2) \rangle
=\langle x',ab_1x_2\rangle+\langle x',ax_1b_2 \rangle\\
&=\langle D(b_2,x_2),x'\widehat{\otimes}ab_1\rangle+\langle
D(b_1,x_1),b_2x'\widehat{\otimes}a\rangle\\
&=\langle D(b_2,x_2),(x'\widehat{\otimes}a)(b_1,x_1)\rangle+\langle
D(b_1,x_1),(b_2,x_2)(x'\widehat{\otimes}a)\rangle\\
&=\langle
(b_1,x_1).(D(b_2,x_2))+(D(b_1,x_1)).(b_2,x_2),x'\widehat{\otimes}a\rangle.
\end{align*}
Thus $D$ is a derivation. Also if
$$(b_\alpha,x_\alpha)~\stackrel{weak^*}{-\hspace{-.2cm}-\hspace{-.2cm}\longrightarrow}
(b,x)\hspace{1cm} in \hspace{0.3cm}  \cal A \oplus_ {\infty}X,$$
then
$x_\alpha~\stackrel{weak^*}{-\hspace{-.2cm}-\hspace{-.2cm}\longrightarrow}
x$ in $X$. Since $X$ is a normal dual $\cal A-$bimodule, then
$ax_\alpha~\stackrel{weak^*}{-\hspace{-.2cm}-\hspace{-.2cm}\longrightarrow}
ax$ in $X$. On the other hand $X$ is reflexive, then
$ax_\alpha~\stackrel{weakly}{-\hspace{-.2cm}-\hspace{-.2cm}\longrightarrow}
ax$ in $X$. Thus $$\langle
D(b_\alpha,x_\alpha),x'\widehat{\otimes}a\rangle=\langle
x',ax_\alpha\rangle \longrightarrow \langle x',ax\rangle=\langle
D(b,x),x'\widehat{\otimes}a\rangle,$$ for every
$x'\widehat{\otimes}a \in {X^*\widehat{\otimes} {\cal A}}.$
Therefore $D$ is $weak^*-$continuous. Connes-amenability of $\cal A
\oplus_ {\infty}X$ implies that $D=\delta_F$ for some $F\in
{(X^*\widehat{\otimes} {\cal A})}^*$. Then for each
$x'\widehat{\otimes}a \in {X^*\widehat{\otimes} {\cal A}}$ and
$(b,x)\in \cal A \oplus_ {\infty}X,$ we have
\begin{align*}\langle x',ax \rangle &= \langle D((b,x)),x'\widehat{\otimes}a
\rangle\\
&=\langle (b,x).F-F.(b,x),x'\widehat{\otimes}a \rangle\\
&=\langle F,(x'\widehat{\otimes}a) .(b,x)-(b,x)(x'\widehat{\otimes}a) \rangle\\
&=\langle F,x'\widehat{\otimes}ab-bx'\widehat{\otimes}a \rangle.\\
\end{align*}
Then $\langle x',ax \rangle=0$ for each $a\in \cal A, x\in X$ and
$x'\in X^*$. We have to  show that $\cal AX=X.$ To this end, we know
that $\cal A \oplus_ {\infty}X$  is Connes-amenable, then it is
unital [5]. Let $(e,x)$ be the unite element of  $\cal A \oplus_
{\infty}X$. It is easy to show that $x=0$ and $ey=y$ for every $y\in
X,$ and the proof is complete.\hfill$\blacksquare~$

\paragraph{\bf Corollary 1.4.}
Let $\cal A$ be a dual Banach algebra and let $X$ be a non-trivial
Banach $\cal A$-bimodule. If $\cal A$ and $X$ are reflexive, then
$\cal A \oplus_ {\infty}X$  is not Connes-amenable.

\paragraph{\bf Corollary 1.5.}
Let $\cal A$ be a non-trivial reflexive dual Banach algebra. Then
the (dual) Banach algebras $\cal A \oplus_ {\infty}\cal A$ and $\cal
A \oplus_ {\infty}\cal A^*$ are not Connes-amenable.
\section{weak Connes-amenability}
Let $\cal A$ be a dual Banach algebra and let $X$ be a normal dual
Banach $\cal A$-bimodule by predual $X_*$ and let $\cal A_*$ be the
predual of $\cal A.$ In lemma 1.1, the module extension Banach
algebra $\cal A \oplus_{\infty} X$  is a dual Banach algebra. We
study the weak Connes-amenability of $\cal A \oplus_{\infty} X$.

\paragraph{\bf Lemma 2.1.}
Let $X$ be a normal, dual Banach $\cal A$-bimodule and
$T:X\longrightarrow X$ be a $weak^*-$continuous $\cal A$-bimodule
morphism. Then $\bar{T}:{\cal A} \oplus_{\infty} X\longrightarrow
{\cal A} \oplus_{\infty} X$, defined by $\ \bar{T}((a,x))=(0,T(x))$
is a $weak^*-$continuous derivation. $\bar{T}$ is inner if and only
if there exists $b\in \cal A$ such that $ba=ab$ for each $a\in \cal
A$ and $T(x)=xb-bx$ for all $x\in X$.
\paragraph{\bf Proof.} Let $(a,x),(b,y)\in {\cal A}
\oplus_{\infty} X$, we have
\[\bar{T}((a,x).(b,y))=\bar{T}((ab,ay+xb))=(0,T(ay+xb))=(0,aT(y))+(0,T(x)b).\]
On the other hand $\bar{T}((a,x)).(b,y)=(0,T(x)).(b,y)=(0,T(x)b)$,
similarly $$(a,x).\bar{T}((b,y))=(a,x).(0,T(x))=(0,aT(y)),$$ and
hence $\bar{T}$ is a derivation. From $weak^*-$continuity of $T$, it
is clear that $\bar{T}$ is $weak^*-$continuous. If $\bar{T}$ is
inner then there exists $\xi=(b,y)\in {\cal A} \oplus_{\infty} X$
such that $\bar{T}((a,x))=(a,x).\xi-\xi.(a,x).$ In particular
$(0,0)=(a,0).\xi-\xi.(a,0)$ and $(0,T(x))=(0,x).\xi-\xi.(0,x).$ Then
$(0,0)=(ab-ba,ay-ya)$ and$(0,T(x))=(0,xb-bx)$ and so there exists
$b\in \cal A$ such that $ba=ab$ for $a\in \cal A$ and $T(x)=xb-bx$
for all $x\in X$. Conversely, if there exists $b\in \cal A$ such
that $ba=ab$ for $a\in \cal A$ and $T(x)=xb-bx$ for all $x\in X$,
then
\[\bar{T}((a,x))=(0,T(x))=(ab-ba,xb-bx)=(a,x).(b,0)-(b,0).(a,x).\] This
shows that $\bar{T}$ is inner, and the proof is
complete.\hfill$\blacksquare~$
\paragraph{\bf Lemma 2.2.}
Let $\cal A$ be a dual Banach algebra and let $X$ be a normal, dual
Banach $\cal A$-bimodule. If $D:{\cal A}\longrightarrow X$ is a
$weak^*-$continuous derivation, then $\bar{D}:({\cal A}
\oplus_{\infty} X)\longrightarrow ({\cal A} \oplus_{\infty} X)$
defined by $\bar{D}((a,x))=(0,D(a))$, is a $weak^*-$continuous
derivation. Furthermore, $\bar{D}$ is inner if and only if $D$ is
inner.
\paragraph{\bf Proof.} It is straightforward to check that $\bar{D}$ is a
$weak^*-$continuous derivation. Now let $\bar{D}$ be inner, then
there exists $\xi=(b,y)\in {\cal A} \oplus_{\infty} X$ such that
$\bar{D}((a,x))=(a,x).\xi-\xi.(a,x).$ In Particular
\[(0,D(a))=\bar{D}((a,0))=(a,0).(b,y)-(b,y).(a,0)=(ab-ba,ay-ya),\]
then $D(a)=ay-ya$ for some $y\in X$ and hence $D$ is inner. The
converse is evident.\hfill$\blacksquare~$
\paragraph{\bf Theorem 2.3.}
Let $\cal A$ be a dual Banach algebra and let $X$ be a normal, dual
Banach $\cal A$-bimodule. Then $\cal {\cal A} \oplus_{\infty} X$ is
weakly Connes-amenable if and only if  the following conditions hold:\\
1. The only $weak^*-$continuous derivations $D:{\cal
A}\longrightarrow \cal A$ for which there is a $weak^*-$continuous
operator $T:X\longrightarrow X$ such that $T(ax)=D(a)x+aT(x)$ and
$T(xa)=xD(a)+T(x)a$ ($a\in \cal A, x\in X$), are the inner derivations.\\
2. $H^1_{w^*}(\cal A,X)=\{0\}$.\\
3. The only $weak^*-$continuous $\cal
A$-bimodule morphism $\Gamma:X\longrightarrow \cal A$ for which
$x\Gamma(y)+\Gamma(x)y=0$ $(x,y\in X),$ is zero.\\
4. For every $weak^*-$continuous $\cal A$-bimodule morphism
$T:X\longrightarrow X$, there exists $b\in \cal A$ for which $ab=ba$
for $a\in \cal A$ and $T(x)=xb-bx$ for $x\in X$.
\paragraph{\bf Proof.} Denote by $\tau_1$
and $\tau_2$ the inclusion mappings from, respectively, $\cal A$ and
$X$ into ${\cal A} \oplus_{\infty} X$, and denote by $\Delta_1$ and
$\Delta_2$ the natural projections from ${\cal A} \oplus_{\infty} X$
onto $\cal A$ and $X$, respectively. Then $\Delta_1$ and  $\Delta_2$
are $\cal A$-bimodule morphisms, so $\tau_1$ and $\tau_2$ are
algebra homomorphisms.To prove the sufficiency we assume that
conditions 1-4 hold. Let $D:{\cal A} \oplus_{\infty}
X\longrightarrow {\cal A} \oplus_{\infty} X$ be a
$weak^*-$continuous derivation. Then $\Delta_1 o D o \tau_1 :{\cal
A}\longrightarrow \cal A$ and $\Delta_2 o D o \tau_1 :{\cal
A}\longrightarrow X$ are $weak^*-$continuous derivations. Now we
show that $\Gamma=\Delta_1 o D o \tau_2 :X\longrightarrow \cal A$ is
trivial. By condition 3 it suffices to show that $\Gamma$ is an
$\cal A$-bimodule morphism satisfying $x\Gamma(y)+\Gamma(x)y=0$
($x,y\in X$). We have
\begin{align*}
0&=D((0,0))=D((0,x).(0,y))=\\
&=D((0,x)).(0,y)+(0,x).D((0,y))\\
&=(0,\Gamma(x)y)+(0,x\Gamma(y)).
\end{align*}
  On the other hand,
\begin{align*}
\Gamma(ax)&=\Delta_1 o D((0,ax))
=\Delta_1oD((a,0).(0,x))\\
&=\Delta_1(D((a,0)).(0,x)+(a,0).D((0,x)))\\
&=\Delta_1((a,0).D((0,x)))
=\Delta_1(aD o \tau_2(x))\\
&=a\Gamma(x).
\end{align*}
Similarly, $\Gamma(xa)=\Gamma(x)a$. Then $\Gamma$ is an $\cal
A$-bimodule morphism such that $x\Gamma(y)+\Gamma(x)y=0.$ Therefore
$\Gamma$ is trivial. Now let $T=\Delta_2 o D o
\tau_2:X\longrightarrow X$ and $D_1=\Delta_1 o D o \tau_1:{\cal
A}\longrightarrow \cal A$. For every $a\in \cal A$ and $x\in X,$ we
have
\begin{align*}(0,T(ax))&=(0,\Delta_2 o
D((0,ax))=D((0,ax))\\
&=D((a,0).(0,x))=D((a,0)).(0,x)+(a,0).D((0,x))\\
&=(0,D_1(a)x)+a(0,T(x))=(0,D_1(a)x+aT(x)) \hspace {1cm}(1).~~~~~
\end{align*}
This means that  $T(ax)=D_1(a)x+aT(x)$. Similarly, for every $a\in
\cal A$ and $x\in X$, we have $$(0,T(xa))=(0,xD_1(a)+T(x)a) \hspace
{1cm}(2).$$
Therefore by condition 1, $D_1=\Delta_1 o D o \tau_1$ is inner.\\
Now suppose that $b\in \cal A$ satisfies $D_1(a)=ab-ba$ for $a\in
\cal A$. Let $T_1:X\longrightarrow X$ be defined by $T_1(x)=xb-bx$
for $x\in X$. Then $T-T_1:X\longrightarrow X$ is a
$weak^*-$continuous $\cal A$-bimodule morphisms. In fact, from (1),
for every $a\in \cal A$ and $x\in X$, we have
\begin{align*} (T-T_1)(ax)&=T(ax)-T_1(ax)\\
&=(D_1(a)x+aT(x))-(axb-bax)\\
&=(ab-ba)x+aT(x)-(axb-bax)\\
&=a(bx-xb)+aT(x)=a(T-T_1)(x).
\end{align*}
 Similarly, $T-T_1$ is a right $\cal
A$-bimodule morphism. From condition 4 there is  $c\in \cal A$ such
that $ac=ca$ for $a\in \cal A$ and $(T-T_1)(x)=xc-cx$ for $x\in X$.
By Lemma 2.1, we know that
\[\overline{T-T_1}:(a,x)\longrightarrow (0,(T-T_1)(x)), {\cal A} \oplus_{\infty}
X\longrightarrow {\cal A} \oplus_{\infty} X\] is an inner
derivation. Since $\Delta_2 o D o \tau_1:{\cal A}\longrightarrow X$
is a $weak^*-$continuous derivation, it is inner by condition 2. By
Lemma 2.2, the mapping
\[\overline{\Delta_2 o D o \tau_1}:(a,x)\longrightarrow (0,\Delta_2 o D o
\tau_1(a)), {\cal A} \oplus_{\infty} X\longrightarrow {\cal A}
\oplus_{\infty} X\] is also inner derivation. Since $\Gamma$ is
trivial, we now have
\begin{align*}D((a,x))&=(D_1(a),\Delta_2 o D o
\tau_1(a)+T(x))\\
&=\overline{\Delta_2 o D o
\tau_1}((a,x))+\overline{(T-T_1)}((a,x))+(D_1(a),T(x)).
\end{align*}
 Since
\[(D_1(a),T_1(x))=(ab-ba,xb-bx)=(a,x).(u,0)-(u,0).(a,x)\] for
$a\in \cal A$ and $x\in X$, it gives an inner derivation from
${\cal A} \oplus_{\infty} X$ into ${\cal A} \oplus_{\infty} X$.
Hence as a sum of three inner derivations, $D$ is inner. Thus
under conditions 1-4,
${\cal A} \oplus_{\infty} X $ is weakly Connes-amenable.\\
Now we prove the necessity. Suppose that ${\cal A} \oplus_{\infty}
X$ is weakly Connes-amenable. Let $D:{\cal A}\longrightarrow \cal A$
be a $weak^*-$continuous derivation with the property given in
condition 1. We define $\bar{D}:{\cal A} \oplus_{\infty}
X\longrightarrow {\cal A} \oplus_{\infty} X$ by
\[\bar{D}((a,x))=(D(a),T(x))\quad(a,x)\in ({\cal A} \oplus_{\infty} X).\]Then $\bar{D}$ is a
$weak^*-$continuous derivation. But $\bar{D}$ is inner, so there
exists $(b,y)\in {\cal A} \oplus_{\infty} X$ such that
$$\bar{D}((a,x))=(a,x).(b,y)-(b,y).(a,x),$$
and then for some $b\in \cal A$, we have
$(D(a),T(x))=(ab-ba,xb-bx)$, thus $D(a)=ab-ba$ , this means that $D$
is inner, and condition 1 holds. Condition 2 follows from Lemma 2.2.
Let now $\Gamma:X\longrightarrow \cal A$ be an arbitrary
$weak^*-$continuous $\cal A$-bimodule morphism for which
$x\Gamma(y)+\Gamma(x)y=0$ ($x,y\in X$). Define $\bar{\Gamma}:{\cal
A} \oplus_{\infty} X\longrightarrow {\cal A} \oplus_{\infty} X$ by
$\bar{\Gamma}((a,x))=(\Gamma(x),0)$ then $\bar{\Gamma}$ is a
$weak^*-$continuous derivation,
 but $\bar{\Gamma}$ is inner, then there exists $\xi=(b,y)\in {\cal A} \oplus_{\infty} X$ such that
 $\bar{\Gamma}((a,x))=(a,x).(b,y)-(b,y).(a,x)$. In particular
  \[(\Gamma(x),0)=\bar{\Gamma}((0,x))=(0,x).(b,y)-(b,y).(0,x)=(0,xb-bx)\]
  and then $\Gamma=0$, and condition 3 holds.
Let $T:X\longrightarrow X$ be a $weak^*-$continuous $\cal
A$-bimodule morphism. $\bar{D}:{\cal A} \oplus_{\infty}
X\longrightarrow {\cal A} \oplus_{\infty} X$ defined by
$\bar{D}((a,x))=(0,T(x))$ is a $weak^*-$continuous derivation, and
condition 4 holds by lemma 2.1.\hfill$\blacksquare~$

 Let $X=\cal A.$ So in condition 4 of above theorem, let
$T=id:\cal A \longrightarrow \cal A$, then we have
\paragraph{\bf Corollary 2.4.} For every non-trivial dual Banach algebra $\cal
A$, we have $H^1_{w^*}({\cal A} \oplus_{\infty} \cal A, {\cal A}
\oplus_{\infty} \cal A)\neq\{o\}.$
\paragraph{\bf Lemma 2.5.} Let $X$ and $Y$ be dual Banach spaces, then
every $weak^*-$continuous linear map from $X$ into $Y$ is bounded.
\paragraph{\bf Proof.} Let $T:X\longrightarrow Y$ be an unbounded linear
map, then there exists a sequence $\{x_n\}$ in $X$ such that
$lim_n\|x_n\|=0$ and $lim_n\|T(x_n)\|=\infty.$  By uniform
boundedness theorem [4],
$T(x_n)~\stackrel{weak^*}{-\hspace{-.2cm}-\hspace{-.2cm}\nrightarrow}
0.$ On the other hand  $weak^*-lim_nx_n=0$, then $T$ is not
$weak^*-$continuous.\hfill$\blacksquare~$

By corollary 2.4 and lemma 2.5, we result
\paragraph{\bf Corollary 2.6.} For every non-trivial dual Banach algebra $\cal
A$, we have $H^1({\cal A} \oplus_{\infty} \cal A, {\cal A}
\oplus_{\infty} \cal A)\neq\{o\}.$

Let $\cal A$ be a dual Banach algebra, and let $X=\cal A$ by module
actions $$a.x=ax, \hspace {1 cm} x.a=0, \hspace {1 cm}(a\in \cal A,
x\in X),$$ we denote $X$ by $\cal A_0,$ then we have the following.
\paragraph{\bf Corollary 2.7.} $\cal A$ is unital and weakly Connes-amenable if and only if ${\cal A} \oplus_{\infty} \cal A_0$ is
weakly Connes-amenable.
\paragraph{\bf Proof.} Let $\cal A$ be weakly Connes-amenable then
the conditions 1 and 2 in theorem 2.3, hold, so if $\cal A$ is
unital then conditions 3 and 4 hold when $X=\cal A_0$. For the
converse let ${\cal A} \oplus_{\infty} \cal A_0$ be weakly
Connes-amenable, then by condition 2, $\cal A$ is weakly
Connes-amenable. The mapping $id:\cal A_0 \longrightarrow \cal A_0$
is a $weak^*-$continuous $\cal A-$bimodule morphism, then by
condition 4 of theorem 2.3, there exists $b\in \cal A$ for which
$ab=ba$ for $a\in \cal A$ and $x=id(x)=x.b-b.x=bx$ for $x\in \cal
A_0$. Thus $b$ is the unite element of $\cal
A.$\hfill$\blacksquare~$

\end{document}